\documentclass{article}
\usepackage{graphicx}
\usepackage{amsthm,amssymb,latexsym,amsmath,xypic}

\newcommand{\R}{\mathbb{R}} \newcommand{\C}{\mathbb{C}}
 
\newcommand{\id}{{\mathbf 1}} 
\newcommand{\dbar}{\overline{\partial}}

\newcommand{\modspc}{{\cal M}_X(E)}
\newcommand{\imodspc}{{\cal M}^*_X(E)}

\newtheorem*{remark}{{\bf Remark}}

\begin{document}

\pagestyle{myheadings}
\markright{Instantons}

\parindent 0mm

\title{Instantons: topological aspects}

\author{Marcos Jardim \\ IMECC - UNICAMP \\
Departamento de Matem\'atica \\ Caixa Postal 6065 \\
13083-970 Campinas-SP, Brazil}

\maketitle

\begin{abstract}

\noindent{\bf Cross-references:} Moduli spaces (62); Instantons in gauge theory (72); 
Mathematical uses of gauge theory (75); Index theorems (149); Donaldson invariants (255);
Characteristic classes (354); Differential geometry (386); Finite dimensional algebras
and quivers (418).
\end{abstract}

\section{Introduction}

Let $X$ be a closed (connected, compact without boundary) smooth
manifold of dimension 4, provided with a Riemannian metric
denoted by $g$. Let $\Omega_X^{p}$ denote space of smooth $p$-forms
on $X$, i.e. the sections of $\wedge^pTX$. The Hodge operator acting
on $p$-forms:
$$ * : \Omega_X^{p} \to \Omega_X^{4-p} $$
satisfies $*^2=(-1)^p$. In particular, $*$ splits $\Omega_X^{2}$
into two sub-spaces $\Omega_X^{2,\pm}$ with eigenvalues $\pm1$:
\begin{equation}\label{dec}
\Omega_X^{2} = \Omega_X^{2,+} \oplus \Omega_X^{2,-} ~~.
\end{equation}
Note also that this decomposition is an orthogonal one, with
respect to the inner product:
$$ \langle \omega_1,\omega_2 \rangle =
\int_X \omega_1\wedge*\omega_2 ~~. $$
A 2-form $\omega$ is said to be self-dual if $*\omega=\omega$ and
it is said to be anti-self-dual if $*\omega=-\omega$. Any 2-form
$\omega$ can be written as the sum
$$ \omega = \omega^+ + \omega^- $$
of its self-dual $\omega^+$ and anti-self-dual $\omega^-$
components.

Now let $E$ be a complex vector bundle over $X$ as above, provided
with a connection $\nabla$, regarded as a $\C$-linear operator
$$ \nabla: \Gamma(E) \to \Gamma(E)\otimes \Omega^1_X $$
satisfying the Leibnitz rule:
$$ \nabla (f\sigma) = f\nabla\sigma + \sigma\otimes df $$
for all $f\in C^{\infty}(X)$ and $\sigma\in\Gamma(E)$. Its curvature
$F_\nabla=\nabla\circ\nabla$ is a 2-form with values  in ${\rm End}(E)$,
i.e. $F_\nabla\in\Gamma({\rm End}(E))\otimes\Omega^2_X$, satisfying the
Bianchi identity $\nabla F_\nabla=0$.

The {\em Yang-Mills equation} is:
\begin{equation}\label{ymeqn}
\nabla * F_\nabla = 0
\end{equation}
It is a 2$^{\rm nd}$-order non-linear equation on the connection
$\nabla$. It amounts to a non-abelian generalization of Maxwell equations,
to which it reduces when $E$ is a line bundle; the four components of
$\nabla$ are interpreted as the electric and magnetic potentials. 

An {\em instanton} on $E$ is a smooth connection $\nabla$ whose curvature
$F_\nabla$ is anti-self-dual as a 2-form, i.e. it satisfies:
\begin{equation}\label{asd}
F_\nabla^+ = 0 ~~{\rm i.e.}~~ *F_\nabla=-F_\nabla ~~.
\end{equation}
The instanton equation is still non-linear (it is linear only if $E$ is a
line bundle), but it is only 1$^{\rm st}$-order on the connection.

Note that if $F_\nabla$ is either self-dual or anti-self-dual as a 2-form,
then the Yang-Mills equation is automatically satisfied:
$$ *F_\nabla=\pm F_\nabla ~~ \Rightarrow ~~
\nabla * F_\nabla = \pm \nabla F_\nabla = 0 $$
by the Bianchi identity. In other words, instantons are particular solutions of
the Yang-Mills equation. Furthermore, while the Yang-Mills equation (\ref{ymeqn})
makes sense over any Riemannian manifold, the instanton equation (\ref{asd}) is
well-defined only in dimension 4.

A gauge transformation is a bundle automorphism $g: E\to E$ covering the identity.
The set of all gauge transformations of a given bundle $E\to X$ form a
group through composition, called the {\em gauge group} and denoted by ${\cal G}(E)$.
The gauge group acts on the set of all smooth connections on $E$ by conjugation:
$$ g\cdot\nabla = g^{-1}\nabla g ~~. $$ 
It is then easy to see that (\ref{asd}) is gauge invariant condition, since
$F_{g\cdot\nabla}=g^{-1}F_\nabla g$. The anti-self-duality equation (\ref{asd}) is also
{\em conformally invariant}: a conformal change in the metric does not change
the decomposition (\ref{dec}), so it preserves self-dual and anti-self-dual 2-forms.

The {\rm topological charge} $k$ of the instanton $\nabla$ is defined by the
integral:
\begin{eqnarray}
\label{k} k &=& - \frac{1}{8\pi^2} \int_X {\rm tr}(F_\nabla\wedge F_\nabla) ~~\\
\nonumber &=& c_2(E) - \frac{1}{2}c_1(E)^2
\end{eqnarray}
where the second equality follows from Chern-Weil theory.

If $X$ is a smooth, non-compact, complete Riemannian manifold, an instanton
on $X$ is an anti-self-dual connection for which the integral (\ref{k}) converges.
Note that in this case, $k$ as above {\em need not be an integer}; however
it is expected to always be quantized, i.e. always a multiple of some fixed
(rational) number which depends only on the base manifold $X$.

\paragraph{Summary.}
This note is organized as follows. After revisiting the variational approach to the 
anti-self-duality equation (\ref{asd}), we study instantons over the simplest
possible Riemannian 4-manifold, $\R^4$ with the flat Euclidean metric. We present 
't Hooft's explicit solutions (Section \ref{euc}), the ADHM construction
(Section \ref{ADHM}) and its dimensional reductions to $\R^3$, $\R^2$ and $\R$
(Section \ref{dimred}). We conclude in Section \ref{moduli} by explaining the
construction of the central object of study in gauge theory, the instanton moduli spaces.

%%%%%%%%%%%%

\section{Variational aspects of Yang-Mills equation.}
Given a fixed smooth vector bundle $E\to X$, let ${\cal A}(E)$ be the set of
all (smooth) connections on $E$. The {\em Yang-Mills functional} is defined by
\begin{equation}\label{ymf} {\rm YM}: {\cal A}(E) \to \R \end{equation}
$$ {\rm YM}(\nabla) = ||F_\nabla||^2_{L^2} = \int_M {\rm tr}(F_\nabla\wedge*F_\nabla) $$
The Euler-Lagrange equation for this functional is exactly the
Yang-Mills equation (\ref{ymeqn}). In particular, self-dual and anti-self-dual
connections yield critical points of the Yang-Mills functional.

Splitting the curvature into its self-dual and anti-self-dual parts, we have
$$ {\rm YM}(\nabla) = ||F_\nabla^+||^2_{L^2} + ||F_\nabla^-||^2_{L^2} $$
It is then easy to see that every anti-self-dual connection $\nabla$ is an absolute
minimum for the Yang-Mills functional, and that ${\rm YM}(\nabla)$ coincides with the
topological charge (\ref{k}) of the instanton $\nabla$ times $8\pi^2$.

One can construct, for various 4-manifolds but most interestingly for $X=S^4$,
solutions of the Yang-Mills equations which are neither self-dual nor anti-self-dual.
Such solutions do not minimize (\ref{ymf}). Indeed, at least for gauge group $SU(2)$
or $SU(3)$, it can be shown that there are no other local minima: any critical point which
is neither self-dual nor anti-self-dual is unstable and must be a ``saddle point"
\cite{BL}. 

%%%%%%%%%%%%

\section{Instantons on Euclidean space}\label{euc}

Let $X=\R^4$ with the flat Euclidean metric, and consider a
hermitian vector bundle $E\to\R^4$. Any connection $\nabla$ on $E$
is of the form $d+A$, where $d$ denotes the usual de Rham operator and
$A\in\Gamma({\rm End}(E))\otimes\Omega^1_{\R^4}$ is a 1-form with values
in the endomorphisms of $E$; this can be written
as follows:
$$ A = \sum_{k=1}^4 A_k dx^k ~~,~~ A_k:\R^4\to\mathfrak{u}(r) ~~. $$
In the Euclidean coordinates $x_1,x_2,x_3,x_4$, the anti-self-duality equation
(\ref{asd}) is given by:
$$ F_{12} = F_{34} ~~~~,~~~~ F_{13} = -F_{24} ~~~~,~~~~ F_{14} = F_{23} $$
where
$$ F_{ij} = \frac{\partial A_j}{\partial x_i} - \frac{\partial A_i}{\partial x_j}
+ [A_i,A_j] ~~.$$

The simplest explicit solution is the charge 1 $SU(2)$-instanton on $\R^4$.
The connection 1-form is given by:
\begin{equation}\label{basic}
A_0 = \frac{1}{1+|x|^2} \cdot {\rm Im}(q d\overline{q})
\end{equation}
where $q$ is the quaternion $q=x_1+x_2{\mathbf i}+x_3{\mathbf j}+x_4{\mathbf k}$,
while ${\rm Im}$ denotes the imaginary part of the product quaternion; we are
regarding ${\mathbf i},{\mathbf j},{\mathbf k}$ as a basis of the Lie algebra
$\mathfrak{su}(2)$; from this, one can compute the curvature:
\begin{equation}
F_{A_0} = \left( \frac{1}{1+|x|^2} \right)^2 \cdot {\rm Im}(dq\wedge d\overline{q})
\end{equation}
We see that the action density function
$$ |F_{A_0}|^2 = \left( \frac{1}{1+|x|^2} \right)^2 $$
has a bell-shaped profile centered at the origin and decaying like $r^{-4}$.

Let $t_{\lambda,y}:\R^4 \to \R^4$ be the isometry given by the composition of
a translation by $y\in\R^4$ with a homothety by $\lambda\in\R^+$. The pullback
connection $t_{\lambda,y}^*A_0$ is still anti-self-dual; more explicitly:
$$ A_{\lambda,y} = t_{\lambda,y}^*A_0 = 
\frac{\lambda^2}{\lambda^2+|x-y|^2} \cdot {\rm Im}(q d\overline{q}) $$
$$ {\rm and}~~ F_{A_{\lambda,y}} = 
\left( \frac{\lambda^2}{\lambda^2+|x-y|^2} \right)^2
\cdot {\rm Im}(dq\wedge d\overline{q}) ~~. $$
Note that the action density function $|F_A|^2$ has again a bell-shaped
profile centered at $y$ and decaying like $r^{-4}$; the parameter
$\lambda$ measures the concentration of the energy density function, and
can be interpreted as the ``size" of the instanton $A_{\lambda,y}$.

Instantons of topological charge $k$ can be obtained by ``superimposing"
$k$ basic instantons, via the so-called {\em 't Hooft ansatz}. Consider
the function $\rho:\R^4\to\R$ given by:
$$ \rho(x) = 1 + \sum_{j=1}^k \frac{\lambda_j^2}{(x-y_j)^2} ~~, $$
where $\lambda_j\in\R$ and $y_j\in\R^4$.
Then the connection 1-form $A=A_\mu dx_\mu$ with coefficients
\begin{equation}\label{mi}
A_\mu = i \sum_{\nu=1}^4 \overline{\sigma}_{\mu\nu}
\frac{\partial}{\partial x_\nu}\ln(\rho(x))
\end{equation}
is anti-self dual; here, $\overline{\sigma}_{\mu\nu}$ are the matrices
given by ($\mu,\nu=1,2,3$):
$$ \overline{\sigma}_{\mu\nu} = \frac{1}{4i}[\sigma_\mu,\sigma_\nu] 
~~ \overline{\sigma}_{\mu4} = \frac{1}{2} \sigma_\mu $$
where $\sigma_\mu$ are the Pauli matrices.

The connection (\ref{mi}) correspond to $k$ instantons centered at points
$y_i$ with size $\lambda_i$. The basic instanton (\ref{basic}) is exactly
(modulo gauge transformation) what one obtains from (\ref{mi}) for the case
$k=1$. The 't Hooft instantons form a $5k$ parameter family of anti-self-dual
connections.

$SU(2)$-instantons are also the building blocks for instantons with general
structure group \cite{BCGW}. Let $G$ be a compact semi-simple Lie group, with
Lie algebra $\mathfrak{g}$. Let $\phi:\mathfrak{su}(2)\to\mathfrak{g}$ be any
injective Lie algebra homomorphism. If $A$ is an anti-self-dual $SU(2)$-connection
1-form, then it is easy to see that $\phi(A)$ is an anti-self-dual $G$-connection
1-form. Using (\ref{mi}) as an example, we have that:
\begin{equation}\label{G-basic}
A = i \sum_{\mu,\nu} \phi(\overline{\sigma}_{\mu\nu})
\frac{\partial}{\partial x_\nu}\ln(\rho(x)) dx_\mu
\end{equation}
is a $G$-instanton on $\R^4$.

While this guarantees the existence of $G$-instantons on $\R^4$, note that the
instanton (\ref{G-basic}) might be reducible (e.g. $\phi$ can simply be the
obvious inclusion of $\mathfrak{su}(2)$ into $\mathfrak{su}(n)$ for any $n$)
and that its charge depends on the choice of representation $\phi$.
Furthermore, it is not clear whether every $G$-instanton can be obtained in
this way, as the inclusion of a $SU(2)$-instanton through some representation
$\phi:\mathfrak{su}(2)\to\mathfrak{g}$.

%%%%%%%%%%%%%%%

\section{The ADHM construction}\label{ADHM}

All $SU(r)$-instantons on $\R^4$ can be obtained through a remarkable construction
due to Atiyah, Drinfeld, Hitchin and Manin. It starts by considering hermitian
vector spaces $V$ and $W$ of dimension $c$ and $r$, respectively, and the
following data:
$$ B_1,B_2 \in {\rm End}(V) ~~ , ~~ i\in{\rm Hom}(W,V)
~~ , ~~ j\in{\rm Hom}(V,W) ~~, $$
so-called ADHM data. Assume moreover that $(B_1,B_2,i,j)$ satisfy the
{\em ADHM equations}:
\begin{eqnarray}
\label{adhm1} [ B_1 , B_2 ] + ij & = & 0  \\
\label{adhm2} [ B_1 , B_1^\dagger ] + [ B_2 , B_2^\dagger ] + ii^\dagger -
j^\dagger j & = & 0
\end{eqnarray}

Now consider the following maps
$$ \alpha:V\times\R^4 \to (V\oplus V\oplus W)\times\R^4 $$
$$ \beta:(V\oplus V\oplus W)\times\R^4 \to V\times\R^4 $$
given as follows ($\id$ denotes the appropriate identity matrix):
\begin{equation} \label{alpha}
\alpha(z_1,z_2) = \left( \begin{array}{c}
B_{1} + z_1\id \\ B_{2} + z_2\id \\ j
\end{array} \right) \end{equation}
\begin{equation} \label{beta}
\beta(z_1,z_2) = \left( \begin{array}{ccc}
-B_{2} - z_2\id \ \ & \ \ B_{1}+ z_1\id \ \ & \ \ i
\end{array} \right) \end{equation}
where $z_1=x_1+ix_2$ and $z_2=x_3+ix_4$ are complex coordinates
on $\R^4$. The maps (\ref{alpha}) and (\ref{beta}) should be understood as a
family of linear maps parameterized by points in $\R^4$.

A straightforward calculation shows that the ADHM equation
(\ref{adhm1}) imply that $\beta\alpha=0$ for every $(z_1,z_2)\in\R^4$.
Therefore the quotient $E=\ker\beta/{\rm im}\alpha=\ker\beta\cap\ker\alpha^\dagger$
forms a complex vector bundle over $\R^4$ or rank $r$ whenever
$(B_1,B_2,i,j)$ is such that $\alpha$ is injective and $\beta$ is
surjective for every $(z_1,z_2)\in\R^4$.

To define a connection on $E$, note that $E$ can be regarded as a sub-bundle
of the trivial bundle $(V\oplus V\oplus W)\times\R^4$. So let
$\iota:E\to(V\oplus V\oplus W)\times\R^4$ be the inclusion, and let
$P:(V\oplus V\oplus W)\times\R^4\to E$ be the orthogonal projection onto $E$.
We can then define a connection $\nabla$ on $E$ through the {\em projection formula}
$$ \nabla s = P \underline{d}\iota(s) $$
where $\underline{d}$ denotes the trivial connection on the trivial bundle
$(V\oplus V\oplus W)\times\R^4$.

To see that this connection is anti-self-dual, note that projection $P$ can be written
as follows
$$ P = \id - {\cal D}^\dagger \Xi^{-1} {\cal D} $$
where
$$ {\cal D} : (V\oplus V\oplus W)\times\R^4 \to (V\oplus V)\times\R^4 $$
$$ {\cal D} =
\left( \begin{array}{c} \beta \\ \alpha^\dagger \end{array} \right) $$
and $\Xi = {\cal D}{\cal D}^\dagger$. Note that $\cal D$ is surjective, so that
$\Xi$ is indeed invertible. Moreover, it also follows from (\ref{adhm2}) that
$\beta\beta^\dagger=\alpha^\dagger\alpha$, so that
$\Xi^{-1}=(\beta\beta^\dagger)^{-1}\id$.

The curvature $F_\nabla$ is given by:
\begin{eqnarray*}
F_{\nabla} & = &
P\left(\underline{d}( \id - {\cal D}_{\rm I}^\dagger\Xi^{-1}{\cal D})\underline{d}\right) = 
P \left(\underline{d}{\cal D}^\dagger\Xi^{-1} (\underline{d}{\cal D})\right) = \\
& = & P \left( (\underline{d}{\cal D}^\dagger) \Xi^{-1} (\underline{d}{\cal D})  +
{\cal D}^\dagger \underline{d}(\Xi^{-1} (\underline{d}{\cal D})  \right) =
% \\ & = & P \left( 
(\underline{d}{\cal D}^\dagger)\Xi^{-1} (\underline{d}{\cal D}) 
% \right)
\end{eqnarray*}
for $P \left({\cal D}^\dagger \underline{d}(\Xi^{-1} (\underline{d}{\cal D}))\right) = 0$
on $E=\ker{\cal D}$. Since $\Xi^{-1}$ is diagonal, we conclude that $F_\nabla$ is 
proportional to $d{\cal D}^\dagger\wedge d{\cal D}$, as a 2-form.

It is then a straightforward calculation to show that each entry of 
$d{\cal D}^\dagger\wedge d{\cal D}$ belongs to
$\Omega^{2,-}$.

The extraordinary accomplishment of Atiyah, Drinfeld, Hitchin and Manin was 
to show that every instanton, up to gauge equivalence, can be obtained in this
way; see e.g. \cite{DK}.
For instance, the basic $SU(2)$-instanton (\ref{basic}) is associated with
the following data ($c=1,r=2)$:
$$ B_1,B_2=0 ~~,~~
i=\left(\begin{array}{c}1\\0\end{array}\right) ~~,~~
j=\left(0 ~~ 1\right) $$

\begin{remark}\rm
The ADHM data $(B_1,B_2,i,j)$ is said to be {\em stable} if $\beta$ is
surjective for every $(z_1,z_2)\in\R^4$, and it is said to be {\em costable} if
$\alpha$ is injective for every $(z_1,z_2)\in\R^4$. $(B_1,B_2,i,j)$ is {\em regular}
if it is both stable and costable. The quotient:
$$ \left. \left\{ \rm{regular ~ solutions ~ of} ~
(\ref{adhm1}) ~ \rm{and} ~ (\ref{adhm2}) \right\} 
\right/ U(V) $$
coincides with the moduli space of instantons of rank $r=\dim W$
and charge $c=\dim V$ on $\R^4$ (see below). It is also an example of
a quiver variety (see the article "Finite dimensional algebras and quivers" by
Alistair Savage), associated to the quiver consisting of two vertices $V$ and
$W$, two loop-edges on the vertex $V$ and two edges linking $V$ to $W$, one in
each direction.
\end{remark}

%%%%%%%%%%%%%

\section{Dimensional reductions of the anti-self-dual Yang-Mills equation.}\label{dimred}

As pointed out above, a connection on an hermitian vector bundle $E\to\R^4$ of rank
$r$ can be regarded as 1-form
$$ A = \sum_{k=1}^4 A_k(x_1,\cdots,x_4) dx^k ~~,~~ A_k:\R^4\to\mathfrak{u}(r) ~~. $$

Assuming that the connection components $A_k$ are invariant under translation
in one direction, say $x_4$, we can think of
$$ \underline{A}=\sum_{k=1}^3 A_k(x_1,x_2,x_3) dx^k $$
as a connection on a hermitian vector bundle over $\R^3$, with the fourth
component $\phi=A_4$ being regarded as a bundle endomorphism
$\phi:E\to E$, called a Higgs field. In this way, the anti-self-duality equation
(\ref{asd}) reduces to the so-called {\em Bogomolny (or monopole) equation}:
\begin{equation} \label{bogomolny}
F_{\underline{A}} = \ast d\phi
\end{equation}
where $\ast$ is the Euclidean Hodge star in dimension 3.

Now assume that the connection components $A_k$ are invariant under translation
in two directions, say $x_3$ and $x_4$. Consider
$$ \underline{A}=\sum_{k=1}^2 A_k(x_1,x_2) dx^k $$
as a connection on a hermitian vector bundle over $\R^2$, with the third and
fourth components combined into a complex bundle endomorphism:
$$ \Phi=(A_3+i\cdot A_4)(dx_1-i\cdot dx_2) $$
taking values on 1-forms.
The anti-self-duality equation (\ref{asd}) is then reduced to the so-called
{\em Hitchin's equations}:
\begin{equation} \label{hitchin} \left\{ \begin{array}{l}
F_{\underline{A}} = [\Phi,\Phi^*] \\ \dbar_{\underline{A}} \Phi = 0
\end{array} \right. \end{equation}
Conformal invariance of the anti-self-duality equation means that
Hitchin's equations are well-defined over any Riemann surface.

Finally, assume that the connection components $A_k$ are invariant under
translation in three directions, say $x_2,x_3$ and $x_4$. After gauging away
the first component $A_1$, the anti-self-duality equations (\ref{asd}) reduce
to the so-called Nahm's equations:
\begin{equation} \label{nahm eqn}
\frac{d T_k}{dx_1} + \frac{1}{2}\sum_{j,l} \epsilon_{kjl}[T_j,T_l] = 0 ~,~~ j,k,l=\{2,3,4\}
\end{equation}
where each $T_k$ is regarded as a map $\R\to\mathfrak{u}(r)$.

Those interested in monopoles and Nahm's equations are referred to the survey \cite{Mu}
and the references therein. The best source for Hitchin's equations still are Hitchin's
original papers \cite{H}. A beautiful duality know as {\em Nahm transform} relates the
various reductions of the anti-self-duality equation to periodic instantons; see the
survey article \cite{J}.

It is also worth mentioning the book by Mason \& Woodhouse \cite{MW},
where other interesting dimensional reductions of the anti-self-duality equation are
discussed, providing a deep relation between instantons and the general theory
of integrable systems.

%%%%%%%%%%%%

\section{The instanton moduli space}\label{moduli}

Now fix a rank $r$ complex vector bundle $E$ over a 4-dimensional Riemannian manifold
$X$. Observe that the difference between any two connections is a linear operator:
$$ (\nabla-\nabla')(f\sigma) = f\nabla\sigma + \sigma\cdot df -
f\nabla'\sigma - \sigma\cdot df = f(\nabla-\nabla')\sigma ~~. $$
In other words, any two connections on $E$ differ by an endomorphism valued 1-form. 
Therefore, the set of all smooth connections on $E$, denoted by ${\cal A}(E)$,
has the structure of an affine space over $\Gamma({\rm End}(E))\otimes\Omega^1_M$.

The gauge group ${\cal G}(E)$ acts on ${\cal A}(E)$ via conjugation:
$$ g\cdot \nabla := g^{-1}\nabla g ~~.$$
We can form the quotient set ${\cal B}(E)={\cal A}(E)/{\cal G}(E)$,
which is the set of gauge equivalence classes of connections on $E$.

The set of gauge equivalence classes of anti-self-dual connections on $E$ is a
subset of ${\cal B}(E)$, and it is called  the {\em moduli space of instantons}
on $E\to X$. The subset of $\modspc$ consisting of irreducible anti-self-dual
connections is denoted $\imodspc$.

Since the choice of a particular vector bundle within its topological class is immaterial,
these sets are usually labeled by the topological invariants (Chern or Pontrjagyn classes)
of the bundle $E$. For instance, ${\cal M}(r,k)$ denotes the moduli space of instantons on
a rank $r$ complex vector bundle $E\to X$ with $c_1(E)=0$ and $c_2(E)=k>0$. 

It turns out that $\modspc$ can be given the structure of a Hausdorff topological space.
In general, $\modspc$ will be singular as a differentiable manifold, but $\imodspc$ can 
always be given the structure of a smooth Riemannian manifold.

We start by explaining the notion of a $L^2_p$ vector bundle. Recall that
$L^2_p(\R^n)$ denotes the completion of the space of smooth functions
$f:\R^n\to \C$ with respect to the norm:
$$ || f ||_{L^2_p}^2 = \int_X \left( |f|^2 + |df|^2 + \cdots |d^{(p)}f|^2 \right) ~~. $$
In dimension $n=4$ and for $p>2$, by virtue of the Sobolev embedding theorem, $L^2_p$
consists of continuous functions, i.e. $L^2_p(\R^n)\subset C^0(\R^n)$. So we define the
notion of a $L^2_p$ vector bundle as a topological vector bundle whose transition functions
are in $L^2_p$, where $p>2$.

Now fixed a $L^2_p$ vector bundle $E$ over $X$, we can consider the metric
space ${\cal A}_p(E)$ of all connections on $E$ which can be represented
locally on an open subset $U\subset X$ as a $L^2_{p}(U)$ 1-form. In this
topology, the subset of irreducible connections ${\cal A}^*_p(E)$ becomes
an open dense subset of ${\cal A}_p(E)$.
Since any topological vector bundle admits a compatible smooth structure, we
may regard $L^2_{p}$ connections as those that differ from a smooth connection
by a $L^2_{p}$ 1-form. In other words, ${\cal A}_p(E)$ becomes an affine space
modeled over the Hilbert space of $L^2_{p}$ 1-forms with values in the endomorphisms
of $E$. The curvature of a connection in ${\cal A}_p(E)$ then becomes a $L^2_{p-1}$
2-form with values in the endomorphism bundle ${\rm End}(E)$.

Moreover, let ${\cal G}_{p+1}(E)$ be defined as the topological group of all
$L^2_{p+1}$ bundle automorphisms. By virtue of the Sobolev multiplication theorem,
${\cal G}_{p+1}(E)$ has the structure of an infinite dimensional Lie group modeled
on a Hilbert space; its Lie algebra is the space of $L^2_{p+1}$ sections of ${\rm End}(E)$.
  
The Sobolev multiplication theorem is once again invoked to guarantee that the
action ${\cal G}_{p+1}(E)\times{\cal A}_p(E)\to{\cal A}_p(E)$ is a smooth map
of Hilbert manifolds. The quotient space ${\cal B}_p(E)={\cal A}_p(E)/{\cal G}_{p+1}(E)$
inherits a topological structure; it is a metric (hence Hausdorff) topological space.
Therefore, the subspace $\modspc$ of ${\cal B}_p(E)$ is also a Hausdorff topological space;
moreover, one can show that the topology of $\modspc$ does not depend on $p$.

The quotient space ${\cal B}_p(E)$ fails to be a Hilbert manifold because in general
the action of ${\cal G}_{p+1}(E)$ on ${\cal A}_p(E)$ is not free. Indeed, if $A$ is a connection on a rank $r$ complex vector bundle $E$ over a connected base manifold $X$,
which is associated with a principal $G$-bundle. Then the isotropy group of $A$ within
the gauge group:
$$ \Gamma_A = \{ g\in{\cal G}_{p+1}(E) ~|~ g(A)=A \} $$
is isomorphic to the centralizer of the holonomy group of $A$ within $G$.

This means that the subspace of irreducible connections ${\cal A}_p^*(E)$ can be
equivalently defined as the open dense subset of ${\cal A}_p(E)$ consisting
of those connections whose isotropy group is minimal, that is:
$$ {\cal A}_p^*(E) = \{ A\in{\cal A}_p(E) ~|~ \Gamma_A={\rm center}(G) \} ~~. $$
Now ${\cal G}_{p+1}(E)$ acts with constant isotropy on ${\cal A}_p^*(E)$,
hence the quotient ${\cal B}_p^*(E)={\cal A}_p^*(E)/{\cal G}_{p+1}(E)$ acquires
the structure of a smooth Hilbert manifold.

\begin{remark}\rm
The analysis of neighborhoods of points in ${\cal B}_p(E)\setminus{\cal B}_p^*(E)$
is very relevant for applications of the instanton moduli spaces to differential topology. 
The simplest situation occurs when $A$ is an $SU(2)$-connection on a rank 2 complex
vector bundle $E$ which reduces to a pair of $U(1)$ and such $[A]$ occurs as an
isolated point in ${\cal B}_p(E)\setminus{\cal B}_p^*(E)$. Then a neighborhood
of $[A]$ in ${\cal B}_p(E)$ looks like a cone on an infinite dimensional complex
projective space.
\end{remark}

Alternatively, the instanton moduli space $\modspc$ can also be described by
first taking the subset of all anti-self-dual connections and then taking the
quotient under the action of the gauge group. More precisely, consider the map:
\begin{equation}\label{map+}
\rho: {\cal A}_p(E) \to L^2_p({\rm End}(E)\otimes\Omega^{2,+}_X)
\end{equation}
$$ \rho(A) = F_A^+ $$
Thus $\rho^{-1}(0)$ is exactly the set of all anti-self-dual connections. It is
${\cal G}_{p+1}(E)$-invariant, so we can take the quotient to get:
$$ \modspc = \rho^{-1}(0)/{\cal G}_{p+1}(E) ~~. $$

It follows that the subspace $\imodspc={\cal B}_p^*(E)\cap\modspc$ has the structure
of a smooth Hilbert manifold. Index theory comes into play to show that $\imodspc$
is finite-dimensional. Recall that if $D$ is an elliptic operator on a vector
bundle over a compact manifold, then $D$ is Fredholm (i.e. $\ker D$ and ${\rm coker}~D$
are finite dimensional) and its index
$$ {\rm ind}~D = \dim\ker D - \dim{\rm coker}~D $$
can be computed in terms of topological invariants, as prescribed by the Atiyah-Singer
index theorem. The goal here is to identify the tangent space of $\imodspc$ with
the kernel of an elliptic operator.

It is clear that for each $A\in{\cal A}_p(E)$, the tangent space $T_A{\cal A}_p(E)$
is just $L^2_p({\rm End}(E)\otimes\Omega^1_X)$. We define the pairing:
\begin{equation}\label{l2pairing}
\langle a , b \rangle = \int_X a \wedge *b
\end{equation}
and it is easy to see that this pairing defines a Riemannian metric
(so-called {\em $L^2$-metric}) on ${\cal A}_p(E)$.

The derivative of the map $\rho$ in (\ref{map+}) at the point $A$ is given by:
$$ d_A^+ ~:~ L^2_{p}({\rm End}(E)\otimes\Omega^1_X) \to 
L^2_{p-1}({\rm End}(E)\otimes\Omega^2_X) $$
$$ a \mapsto (d_A a)^+ ~~, $$
so that for each $A\in\rho^{-1}(0)$ we have:
$$ T_A\rho^{-1}(0) = \left\{ a \in L^2_p({\rm End}(E))\otimes\Omega^1_X ~|~
d_A^+a=0 \right\} ~~. $$

Now for a gauge equivalence class $[A]\in{\cal B}_p^*(E)$, the tangent space
$T_{[A]}{\cal B}_p^*(E)$ consists of those 1-forms which are orthogonal to
the fibers of the principal ${\cal G}_{p+1}(E)$ bundle
${\cal A}_p^*(E)\to{\cal B}_p^*(E)$. At a point $A\in{\cal A}_p(E)$, the
derivative of the action by some $g\in{\cal G}_{p+1}(E)$ is
$$ -d_A ~:~ L^2_{p+1}({\rm End}(E)) \to L^2_p({\rm End}(E)\otimes\Omega^1_X) ~~. $$

Usual Hodge decomposition gives us that there is an orthogonal decomposition:
$$ L^2_p({\rm End}(E)\otimes\Omega^1_X) =
{\rm im}~d_A \oplus \ker d_A^* ~~,$$
which means that:
$$ T_{[A]}{\cal B}_p^*(E) =
\left\{ a \in L^2_p({\rm End}(E)\otimes\Omega^1_X) ~|~ d_A^*a=0 \right\} ~~. $$
Thus the pairing (\ref{l2pairing}) also defines a Riemannian metric on
${\cal B}_p^*(E)$.

Putting these together, we conclude that the space $T_{[A]}\imodspc$ tangent to
$\imodspc$ at an equivalence class $[A]$ of anti-self-dual connections can be
described as follows:
\begin{equation}
T_{[A]}\imodspc = \left\{ a \in L^2_p({\rm End}(E)\otimes\Omega^1_X) ~|~
d_A^*a=d_A^+a=0 \right\}
\end{equation}
It turns out that the so-called {\em deformation operator}
$\delta_A=d_A^* \oplus d_A$:
$$ \delta_A ~:~ 
L^2_p({\rm End}(E)\otimes\Omega^1_X) \to
L^2_{p+1}({\rm End}(E)) \oplus L^2_{p-1}({\rm End}(E)\otimes\Omega^2_X) $$
is elliptic. Moreover, if $A$ is anti-self-dual then ${\rm coker}~\delta_A$
is empty, so that $T_{[A]}\imodspc=\ker\delta_A$. The dimension of the
tangent space $T_{[A]}\imodspc$ is then simply given by the index of
the deformation operator $\delta_A$. Using the Atiyah-Singer index theorem,
we have for $SU(r)$-bundles with $c_2(E)=k$:
$$ \dim \imodspc = 4rk - (r^2-1)(1-b_1(X)+b_+(X)) ~~. $$
The dimension formula for arbitrary gauge group $G$ can be found at
\cite{AHS}.

For example, the moduli space of $SU(2)$ instantons on $\R^4$ of charge $k$
is a smooth Riemannian manifold of dimension $8k-3$. These parameters are
interpreted as the $5k$ parameters describing the positions and sizes of
$k$ separate instantons, plus $3(k-1)$ parameters describing their relative
$SU(2)$ phases. 

The detailed construction of the instanton moduli spaces can be found at
\cite{DK}. An alternative source is Morgan's lecture notes in \cite{FrMo}.

It is interesting to note that $\imodspc$ inherits many of the geometrical
properties of the original manifold $X$. Most notably, if $X$ is a K\"ahler manifold,
then $\imodspc$ is also K\"ahler; if $X$ is a hyperk\"ahler manifold, then $\imodspc$
is also hyperk\"ahler. One expects that other geometric structures on $X$ can also be
transfered to the instanton moduli spaces $\imodspc$.

%%%%%%%%%%%%


\begin{thebibliography}{99}

\bibitem{AHS}
M. F. Atiyah, N. J. Hitchin \& I. M. Singer.
Self-duality in four-dimensional riemannian geometry.
Proc. Royal Soc. London {\bf 362}, 425-461 (1978).

\bibitem{BCGW}
C. W. Bernard, N. H. Christ, A. H. Guth \& E. J. Weinberg.
Pseudoparticle parameters for arbitrary gauge groups.
Phys. Rev. D {\bf 16}, 2967-2977 (1977).

\bibitem{BL}
J. P. Bourguignon \& H. B. Lawson Jr.
Stability and isolation phenomena for Yang-Mills fields
Commun. Math. Phys. {\bf 79} 189-230 (1981).

\bibitem{DK}
S. K. Donaldson \& P. B. Kronheimer.
Geometry of four-manifolds.
Oxford: Clarendon Press 1990.

\bibitem{FrMo}
R. Friedman \& J. W. Morgan (eds.)
Gauge theory and the topology of four-manifolds.
Providence, RI: American Mathematical Society 1998.

\bibitem{H}
N. Hitchin.
The self-duality equations on a Riemann surface.
Proc. London Math. Soc. {\bf 55}, 59-126 (1987). \newline
Stable bundles and integrable systems.
Duke Math. J. {\bf 54}, 91-114 (1987).

\bibitem{J}
M. Jardim.
A survey on Nahm transform.
J. Geom. Phys. {\bf 52}, 313-327 (2004).

\bibitem{MW}
L. J. Mason \& N. M. J. Woodhouse.
Integrability, self-duality, and twistor theory.
New York, NY: Clarendon Press 1996.

\bibitem{Mu}
M. Murray.
Monopoles.
In: Geometric analysis and applications to quantum field theory, 119-135. 
Progr. Math. {\bf 205}. 
Boston, MA: Birkhäuser Boston 2002. 

\end{thebibliography}
\end{document}